\documentclass{amsart}
\usepackage[doi=false, isbn=false, url=false, style=alphabetic, backend=bibtex]{biblatex}
\AtEveryBibitem{\clearfield{MRCLASS}} 
\AtEveryBibitem{\clearfield{MRNUMBER}} 
\AtEveryBibitem{\clearfield{MRREVIEWER}}  
\usepackage{accents}
\usepackage{graphicx}
\usepackage{amssymb}
\usepackage{amsmath}
\usepackage{mathtools}
\makeatletter
\providecommand*{\twoheadrightarrowfill@}{%
  \arrowfill@\relbar\relbar\twoheadrightarrow
}
\providecommand*{\twoheadleftarrowfill@}{%
  \arrowfill@\twoheadleftarrow\relbar\relbar
}
\providecommand*{\xtwoheadrightarrow}[2][]{%
  \ext@arrow 0579\twoheadrightarrowfill@{#1}{#2}%
}
\providecommand*{\xtwoheadleftarrow}[2][]{%
  \ext@arrow 5097\twoheadleftarrowfill@{#1}{#2}%
}
\makeatother
\usepackage{amsfonts}
\usepackage{amscd}
\usepackage{url}
\usepackage{tikz}
\usepackage{calligra}
\usetikzlibrary{shapes}
\usetikzlibrary{decorations.markings}
\usepackage{tikz-cd}
\usetikzlibrary{arrows}
\usepackage{multirow}
\usepackage[savepos]{zref}
\usepackage{hyperref}
\newtheorem{thm}{Theorem}[section]
\newtheorem{corollary}[thm]{Corollary}
\newtheorem{lemma}[thm]{Lemma}

\theoremstyle{definition}

\theoremstyle{remark}
\newtheorem{remark}[thm]{Remark}

\numberwithin{equation}{section}
\numberwithin{figure}{section}

\usepackage[mathscr]{euscript}
\usepackage[T2A,T1]{fontenc}

  {\begin{description}%
    \setlength{\itemsep}{2.5pt}%
    \setlength{\parskip}{5pt}}%
  {\end{description}}

\newcommand{\form}{\upsilon}

\DeclareMathOperator{\SC}{sc}

\DeclareMathOperator{\Spl}{Spl}

\DeclareMathOperator{\Hom}{Hom}

\DeclareMathOperator{\Gal}{Gal}

\DeclareMathOperator{\Ker}{Ker}

\DeclareMathOperator{\Spec}{Spec}

\DeclareMathOperator{\et}{\acute{e}t}

\newcommand{\From}{\colon}
\newcommand{\inar}{\ar@{^{(}->}}
\newcommand{\onar}{\ar@{->>}}

\newcommand{\defined}[1]{\underline{{#1}}}
\usepackage{accents}
\newlength{\dtildeheight}

\newcommand{\Baer}{\dotplus}

\newcommand{\Weil}{\mathcal{W}}

\makeatletter
\newcommand{\raisemath}[1]{\mathpalette{\raisem@th{#1}}}
\newcommand{\raisem@th}[3]{\raisebox{#1}{$#2#3$}}
\makeatother

\newcommand{\EL}{{\mathsf L}}

\newcommand{\Cat}[1]{ {\mathsf{#1}} }

\newcommand{\alg}[1]{\boldsymbol{\mathrm{#1}}}

\newcommand{\sheaf}[1]{{\mathscr{#1}}}

\newcommand{\mGal}{\widetilde{\Gal}}

\newcommand{\gerb}[1]{\boldsymbol{\Cat{#1}}}

\newcommand{\ZZ}{\mathbb Z}

\newcommand{\CC}{\mathbb C}

\newcommand{\Into}{\hookrightarrow}
\newcommand{\Onto}{\twoheadrightarrow}
\newcommand{\To}{\rightarrow}

\newcommand{\inarrow}{\arrow[hook]}
\newcommand{\onarrow}{\arrow[two heads]}

\makeatletter
\newcommand\@biprod[1]{%
  \vcenter{\hbox{\ooalign{$#1\prod$\cr$#1\coprod$\cr}}}}
\newcommand\biprod{\mathop{\mathpalette\@biprod\relax}\displaylimits}
\makeatother

\DeclareMathAlphabet{\mathcalligra}{T1}{calligra}{m}{n}

\DeclareMathOperator{\Wh}{Whit}
\begin{document}

\title{A comparison of L-groups for covers of split reductive groups}%
\author{Martin H. Weissman}%
\date{\today}

\address{Yale-NUS College, 6 College Ave East, \#B1-01, Singapore 138614}
\email{marty.weissman@yale-nus.edu.sg}%


\begin{abstract}
In one article, the author has defined an L-group associated to a cover of a quasisplit reductive group over a local or global field.  In another article, Wee Teck Gan and Fan Gao define (following an unpublished letter of the author) an L-group associated to a cover of a pinned split reductive group over a local or global field.  In this short note, we give an isomorphism between these L-groups.  In this way, the results and conjectures discussed by Gan and Gao are compatible with those of the author.  Both support the same Langlands-type conjectures for covering groups.
\end{abstract}

\maketitle

\tableofcontents

\section*{Summary of two constructions}

Let $\alg{G}$ be a split reductive group over a local or global field $F$.  Choose a Borel subgroup $\alg{B} = \alg{T} \alg{U}$ containining a split maximal torus $\alg{T}$ in $\alg{G}$.  Let $X = \Hom(\alg{T}, \alg{G}_m)$ be the character lattice, and $Y = \Hom(\alg{G}_m, \alg{T})$ be the cocharacter lattice of $\alg{T}$.  Let $\Phi \subset X$ be the set of roots and $\Delta$ the subset of simple roots.  For each root $\alpha \in \Phi$, let $\alg{U}_\alpha$ be the associated root subgroup.  Let $\Phi^\vee$ and $\Delta^\vee$ be the associated coroots and simple coroots.  The root datum of $\alg{G} \supset \alg{B} \supset \alg{T}$ is
$$\Psi = (X, \Phi, \Delta, Y, \Phi^\vee, \Delta^\vee).$$
Fix a pinning (\'epinglage) of $\alg{G}$ as well -- a system of isomorphisms $x_\alpha \From \alg{G}_a \To \alg{U}_\alpha$ for every root $\alpha$.


The following notions of covering groups and their dual groups match those in \cite{MWBig}.  Let $\alg{\tilde G} = (\alg{G}', n)$ be a degree $n$ cover of $\alg{G}$ over $F$; in particular, $\# \mu_n(F) = n$.  Here $\alg{G}'$ is a central extension of $\alg{G}$ by $\alg{K}_2$ in the sense of \cite{B-D}, and write $(Q, \sheaf{D}, f)$ for the three Brylinski-Deligne invariants of $\alg{G}'$.  Assume that if $n$ is odd, then $Q \From Y \To \ZZ$ takes only even values (this is \cite[Assumption 3.1]{MWBig}).

Let $\tilde G^\vee \supset \tilde B^\vee \supset \tilde T^\vee$ be the dual group of $\alg{\tilde G}$, and let $\tilde Z^\vee$ be the center of $\tilde G^\vee$.  The group $\tilde G^\vee$ is a pinned complex reductive group, associated to the root datum
$$(Y_{Q,n}, \tilde \Phi^\vee, \tilde \Delta^\vee, X_{Q,n}, \tilde \Phi, \tilde \Delta).$$
Here $Y_{Q,n} \subset Y$ is a sublattice containing $nY$.  For each coroot $\alpha^\vee \in \Phi^\vee$, there is an associated positive integer $n_\alpha$ dividing $n$ and a ``modified coroot'' $\tilde \alpha^\vee = n_\alpha \alpha^\vee \in \tilde \Phi^\vee$.  The set $\tilde \Phi^\vee$ consists of the modified coroots, and $\tilde \Delta^\vee$ the modified simple coroots.  Define $Y_{Q,n}^{\SC}$ to be the sublattice of $Y_{Q,n}$ generated by the modified coroots.  Then
$$\tilde T^\vee = \Hom(Y_{Q,n}, \CC^\times) \text{ and } \tilde Z^\vee = \Hom(Y_{Q,n} / Y_{Q,n}^{\SC}, \CC^\times).$$

Let $\bar F / F$ be a separable algebraic closure, and $\Gal_F = \Gal(\bar F / F)$ the absolute Galois group.  Fix an injective character $\epsilon \From \mu_n(F) \hookrightarrow \CC^\times$.  From this data, the constructions of \cite{MWBig} and \cite{GanGao} both yield an L-group of $\alg{\tilde G}$ via a Baer sum of two extensions.  In both papers, an extension
\begin{equation}
\tag{First twist} \tilde Z^\vee \Into E_1 \Onto \Gal_F
\end{equation}
is described in essentially the same way  When $F$ is local, this ``first twist'' $E_1$ is defined via a $\tilde Z^\vee$-valued 2-cocycle on $\Gal_F$.  See \cite[\S 5.2]{GanGao} and \cite[\S 5.4]{MWBig} (in the latter, $E_1$ is denoted $(\tau_Q)_\ast \mGal_F$).  Over global fields, the construction follows from the local construction and Hilbert reciprocity.

Both papers include a ``second twist''.  Gan and Gao \cite[\S 5.2]{GanGao} describe an extension
\begin{equation}
\tag{Second twist}
\tilde Z^\vee \Into E_2 \Onto \Gal_F,
\end{equation}
following an unpublished letter (June, 2012) from the author to Deligne.  In \cite{MWBig}, the second twist is the fundamental group of a gerbe, denoted $\pi_1^{\et}(\gerb{E}_\epsilon(\alg{\tilde G}), \bar s)$.  In this article $\bar s = \Spec(\bar F)$, and so we write $\pi_1^{\et}(\gerb{E}_\epsilon(\alg{\tilde G}), \bar F)$ instead.

Both papers proceed by taking the Baer sum of these two extensions, $E = E_1 \Baer E_2$, to form an extension $\tilde Z^\vee \Into E \Onto \Gal_F$.  The extension $E$ is denoted ${}^\EL \tilde Z$ in \cite[\S 5.4]{MWBig}.  Then, one pushes out the extension $E$ via $\tilde Z^\vee \Into \tilde G^\vee$, to define the L-group
\begin{equation}
\tag{L-group}
\tilde G^\vee \Into {}^\EL \tilde G \Onto \Gal_F.
\end{equation}

The two constructions of the L-group, from \cite{GanGao} and \cite{MWBig} are the same, except for insignificant linguistic differences, and a significant difference between the ``second twists''.  In this short note, by giving an isomorphism,
$$ \pi_1^{\et}(\gerb{E}_\epsilon(\alg{\tilde G}), \bar F) \text{ (described by the author)} \xrightarrow{\sim} E_2 \text{ (described by Gan and Gao) }$$
we will demonstrate that the second twists, and thus the L-groups, of both papers are isomorphic.  Therefore, the work of Gan and Gao in \cite{GanGao} supports the broader conjectures of \cite{MWBig}.

\begin{remark}
Among the ``insignificant linguistic differences,'' we note that Gan and Gao use extensions of $F^\times / F^{\times n}$ (for local fields) or the Weil group $\Weil_F$ rather than $\Gal_F$.  But pulling back via the reciprocity map of class field theory yields extensions of $\Gal_F$ by $\tilde Z^\vee$ as above.
\end{remark}

\section{Computations in the gerbe}

\subsection{Convenient base points}

Let $\gerb{E}_\epsilon(\alg{\tilde G})$ be the gerbe constructed in \cite[\S 3]{MWBig}.  Rather than using the language of \'etale sheaves over $F$, we work with $\bar F$-points and trace through the $\Gal_F$-action.  Let $\hat T = \Hom(Y_{Q,n}, \bar F^\times)$ and $\hat T_{\SC} = \Hom(Y_{Q,n}^{\SC}, \bar F^\times)$.  Let $p \From \hat T \To \hat T_{\SC}$ be the surjective $\Gal_F$-equivariant homomorphism dual to the inclusion $Y_{Q,n}^{\SC} \Into Y_{Q,n}$.  Define 
$$\hat Z = \Ker(p) = \Hom(Y_{Q,n} / Y_{Q,n}^{\SC}, \bar F^\times).$$
The reader is warned not to confuse $\hat T, \hat T_{\SC}, \hat Z$ with $\tilde T^\vee, \tilde T_{\SC}^\vee, \tilde Z^\vee$; the former are nontrivial $\Gal_F$-modules (Homs into $\bar F^\times$) and the latter are trivial $\Gal_F$-modules (Homs into $\CC^\times$ as a trivial $\Gal_F$-module).

Write $\bar D = \sheaf{D}(\bar F)$ and $D = \sheaf{D}(F)$, where we recall $\sheaf{D}$ is the second Brylinski-Deligne invariant of the cover $\alg{\tilde G}$.  We have a $\Gal_F$-equivariant short exact sequence,
$$\bar F^\times \Into \bar D \Onto Y.$$
By Hilbert's Theorem 90, the $\Gal_F$-fixed points give a short exact sequence,
$$F^\times \Into D \Onto Y.$$

Let $\bar D_{Q,n}$ and $\bar D_{Q,n}^{\SC}$ denote the preimages of $Y_{Q,n}$ and $Y_{Q,n}^{\SC}$ in $\bar D$.  These are {\em abelian} groups, fitting into a commutative diagram with exact rows.
$$
\begin{tikzcd}
\bar F^\times \inarrow{r} \arrow{d}{=} & \bar D_{Q,n}^{\SC} \onarrow{r} \inarrow{d} & Y_{Q,n}^{\SC} \inarrow{d} \\
\bar F^\times \inarrow{r} & \bar D_{Q,n} \onarrow{r} & Y_{Q,n} 
\end{tikzcd}
$$
Let $\Spl(\bar D_{Q,n})$ be the $\hat T$-torsor of splittings of $\bar D_{Q,n}$, and similarly let $\Spl(\bar D_{Q,n}^{\SC})$ be the $\hat T_{\SC}$-torsor of splittings of $\bar D_{Q,n}^{\SC}$.

Let $\overline{\Wh}$ denote the $\hat T_{\SC}$-torsor of nondegenerate characters of $\alg{U}(\bar F)$.  An element of $\overline{\Wh}$ is a homomorphism (defined over $\bar F$) from $\alg{U}$ to $\alg{G}_a$ which is nontrivial on every simple root subgroup $\alg{U}_\alpha$.  $\Gal_F$ acts on $\overline{\Wh}$, and the fixed points $\Wh = \overline{\Wh}^{\Gal_F}$ are those homomorphisms from $\alg{U}$ to $\alg{G}_a$ which are defined over $F$.  The $\hat T_{\SC}$-action on $\overline{\Wh}$ is described in \cite[\S 3.3]{MWBig}.

The pinning $\{ x_\alpha : \alpha \in \Phi \}$ of $\alg{G}$ gives an element $\psi \in \Wh$.  Namely, let $\psi$ be the unique nondegenerate character of $\alg{U}$ which satisfies
$$\psi( x_\alpha(1) ) = 1 \text{ for all } \alpha \in \Delta.$$

In \cite[\S 3.3]{MWBig}, we define an surjective homomorphism $\mu \From \hat T_{\SC} \To \hat T_{\SC}$, and a $\Gal_F$-equivariant isomorphism of $\hat T_{\SC}$-torsors,
$$\bar \omega \From \mu_\ast \overline{\Wh} \To \Spl(D_{Q,n}^{\SC}).$$
The isomorphism $\bar \omega$ sends $\psi$ to the unique splitting $s_\psi \in \Spl(D_{Q,n}^{\SC})$ which satisfies
$$s_\psi(\tilde \alpha^\vee) = r_\alpha \cdot [e_{\alpha}]^{n_\alpha}, \text{ with } r_\alpha = (-1)^{ Q(\alpha^\vee) \cdot \frac{n_\alpha(n_\alpha-1)}{2} }.$$

We describe the element $[e_{\alpha}] \in D$ concisely here, based on \cite[\S 11]{B-D} and \cite[\S 2.4]{GanGao}.  Let $F(\!(\form)\!)$ be the field of Laurent series with coefficients in $F$.  The extension $\alg{K}_2 \Into \alg{G}' \Onto \alg{G}$ splits over any unipotent subgroup, and so the pinning homomorphisms $x_\alpha \From F(\!(\form)\!)\To \alg{U}_\alpha(F(\!(\form)\!))$ lift to homomorphisms
$$\tilde x_\alpha \From F(\!(\form)\!) \To \alg{U}_\alpha'(F(\!(\form)\!)).$$
Define, for any $u \in F(\!(\form)\!)^\times$,
$$\tilde n_\alpha(u) = \tilde x_\alpha(u) \tilde x_{-\alpha}(-u^{-1}) \tilde x_\alpha(u).$$  
This yields an element
$$\tilde t_\alpha = \tilde n_\alpha( \form) \cdot \tilde n_\alpha(-1) \in \alg{T}'( F(\!(\form)\!)).$$
Then $t_{\alpha}$ lies over $\alpha^\vee(\form) \in \alg{T}( F(\!(\form)\!))$.  Its pushout via $\alg{K}_2( F(\!(\form)\!)) \xtwoheadrightarrow{\partial} F^\times$ is the element we call $[e_\alpha] \in D$.

\begin{remark}
The element $s_\psi(\tilde \alpha^\vee) = r_\alpha \cdot [e_{\alpha}]^{n_\alpha}$ coincides with what Gan and Gao call $s_{Q^{\SC}}(\tilde \alpha^\vee)$ in \cite[\S 5.2]{GanGao}; the sign $r_\alpha$ arises from the formulae of \cite[\S 11.1.4, 11.1.5]{B-D}.
\end{remark}

Let $j_0 \From \hat T_{\SC} \To \mu_\ast \overline{\Wh}$ be the unique isomorphism of $\hat T_{\SC}$-torsors which sends $1$ to $\psi$ (or rather the image of $\psi$ via $\overline{\Wh} \To \mu_\ast \overline{\Wh}$).  Since $\psi \in \Wh$ is $\Gal_F$-invariant, this isomorphism $j_0$ is also $\Gal_F$-invariant.

Finally, let $s \in \Spl(\bar D_{Q,n})$ be a splitting which restricts to $s_\psi$ on $Y_{Q,n}^{\SC}$.  Such a splitting $s$ exists, since the map $\Spl(\bar D_{Q,n}) \To \Spl(\bar D_{Q,n}^{\SC})$ is surjective (since the map $\hat T \To \hat T_{\SC}$ is surjective).  Note that $s$ is not necessarily $\Gal_F$-invariant (and often cannot be).

Let $h \From \hat T \To \Spl(\hat D_{Q,n})$ be the function given by
$$h(x) = x^n \ast s \text{ for all } x \in \hat T.$$

The triple $\bar z = (\hat T, h, j_0)$ is an $\bar F$-object (i.e., a geometric base point) of the gerbe $\gerb{E}_\epsilon(\alg{\tilde G})$.  Note that the construction of $\bar z$ depends on two choices:  a pinning of $\alg{G}$ (to obtain $\psi \in \Wh$) and a splitting $s$ of $\bar D_{Q,n}$ extending $s_\psi$.  We call such a triple $\bar z$ a \defined{convenient base point} for the gerbe $\gerb{E}_\epsilon(\alg{\tilde G})$.

\subsection{The fundamental group}

For a convenient base point $\bar z$ associated to $s$, we consider the fundamental group
$$\pi_1^{\et}(\gerb{E}_\epsilon(\alg{\tilde G}), \bar z) = \bigsqcup_{\gamma \in \Gal_F} \Hom(\bar z, {}^\gamma \bar z).$$
This fundamental group fits into a short exact sequence
$$\tilde Z^\vee \Into \pi_1^{\et}(\gerb{E}_\epsilon(\alg{\tilde G}), \bar z)  \Onto \Gal_F,$$
where the fibre over $\gamma \in \Gal_F$ is $\Hom(\bar z, {}^\gamma \bar z)$.  Thus to describe the fundamental group, it suffices to describe each fibre (as a $\tilde Z^\vee$-torsor), and the multiplication maps among fibres.

The base point ${}^\gamma \bar z$ is the triple $({}^\gamma \hat T, \gamma \circ h, \gamma \circ j_0)$, where ${}^\gamma \hat T$ is the $\hat T$-torsor with underlying set $\hat T$ and twisted action
$$u \ast_\gamma x = \gamma^{-1}(u) \cdot x.$$

To give an element $f \in \Hom(\bar z, {}^\gamma \bar z)$ is the same as giving an element $\zeta \in \tilde Z^\vee$ and a map of $\hat T$-torsors $f_0 \From \hat T \to {}^\gamma \hat T$ satisfying 
$$(\gamma \circ h) \circ f_0 = h \text{ and } (\gamma \circ j_0) \circ p_\ast f_0 = j_0.$$  
Any such map of $\hat T$-torsors is uniquely determined by the element $\tau \in \hat T$ satisfying $f_0(1) = \tau$.  The two conditions above are equivalent to the two conditions
\begin{equation}
\tau^n = \gamma^{-1} s / s \text{ and } \tau \in \hat Z.
\end{equation}

Thus, to give an element $f \in \Hom(\bar z, {}^\gamma \bar z)$ is the same as giving a pair $(\tau, \zeta) \in \hat T \times \tilde Z^\vee$, where $\tau$ satisfies the two conditions above.  Therefore, in what follows, we write $(\tau, \zeta) \in \Hom(\bar z, {}^\gamma \bar z)$ to indicate that $\tau$ satisfies the two conditions above, and to refer to the corresponding morphism in the gerbe $\gerb{E}_\epsilon(\alg{\tilde G})$ in concrete terms.  

We use $\epsilon \From \mu_n(F) \xrightarrow{\sim} \mu_n(\CC)$ to identify $\hat Z_{[n]}$ with $\tilde Z_{[n]}^\vee$.  Two pairs $(\tau, \zeta)$ and $(\tau', \zeta')$ are identified in $\Hom(\bar z, {}^\gamma \bar z)$ if and only if there exists $\xi \in \hat Z_{[n]}$ such that
$$\tau' = \xi \cdot \tau \text{ and } \zeta' = \epsilon(\xi)^{-1} \cdot \zeta.$$

The structure of $\Hom(\bar z, {}^\gamma \bar z)$ as a $\tilde Z^\vee$-torsor is by scaling the second factor in $(\tau, \zeta) \in \hat T \times \tilde Z^\vee$.  To describe the fundamental group completely, it remains to describe the multiplication maps among fibres.  If $\gamma_1, \gamma_2 \in \Gal_F$, and 
$$(\tau_1, \zeta_1) \in \Hom(\bar z, {}^{\gamma_1} \bar z) \text{ and } (\tau_2, \zeta_2) \in \Hom(\bar z, {}^{\gamma_2} \bar z),$$
then their composition in $\pi_1^{\et}(\gerb{E}_\epsilon(\alg{\tilde G}), \bar z)$ is given by
$$(\tau_1, \zeta_1) \circ (\tau_2, \zeta_2) = (\gamma_2^{-1}(\tau_1) \cdot \tau_2, \zeta_1 \zeta_2).$$
Observe that
$$(\gamma_2^{-1}(\tau_1) \tau_2)^n = \gamma_2^{-1} \left( \gamma_1^{-1} s / s \right) \cdot \left( \gamma_2^{-1} s / s \right) =(\gamma_1 \gamma_2)^{-1} s / s.$$
Therefore $(\gamma_2^{-1}(\tau_1) \cdot \tau_2, \zeta_1 \zeta_2) \in \Hom(\bar z, {}^{\gamma_1 \gamma_2} \bar z)$ as required.

\section{Comparison to the second twist}

\subsection{The second twist}

The construction of the second twist in \cite{GanGao} does not rely on gerbes at all, at the expense of some generality; it seems difficult to extend the construction there to nonsplit groups.  But for split groups, the construction of \cite{GanGao} offers significant simplifications over \cite{MWBig}.  The starting point in \cite{GanGao} is the same short exact sequence of abelian groups as in the previous section,
$$F^\times \Into D_{Q,n} \Onto Y_{Q,n}.$$
And as before, we utilize the splitting $s_\psi \From Y_{Q,n}^{\SC} \Into D_{Q,n}^{\SC}$.  Taking the quotient by $s_\psi(Y_{Q,n}^{\SC})$, we obtain a short exact sequence
$$F^\times \Into \frac{D_{Q,n}}{s_\psi(Y_{Q,n}^{\SC})} \Onto \frac{Y_{Q,n}}{Y_{Q,n}^{\SC}}.$$
Apply $\Hom(\bullet, \CC^\times)$ (and note $\CC^\times$ is divisible) to obtain a short exact sequence,
$$\tilde Z^\vee \Into \Hom \left( \frac{D_{Q,n}}{s_\psi(Y_{Q,n}^{\SC})}, \CC^\times \right) \Onto \Hom(F^\times, \CC^\times).$$
Define a homomorphism $\Gal_F \To \Hom(F^\times, \CC^\times)$ by the Artin symbol,
$$\gamma \mapsto \left( u \mapsto \epsilon \left( \frac{ \gamma^{-1}( \sqrt[n]{u}) }{ \sqrt[n]{u} } \right) \right).$$
Pulling back the previous short exact sequence by this homomorphism yields a short exact sequence
$$\tilde Z^\vee \Into E_2 \Onto \Gal_F.$$
This $E_2$ is the second twist described in \cite{GanGao}.

\begin{remark}
There is an insignificant difference here -- at the last step, over a local field $F$, Gan and Gao pull back to $F^\times / F^{\times n}$ via the Hilbert symbol whereas we pull further back to $\Gal_F$ via the Artin symbol.
\end{remark}

Write $E_{2,\gamma}$ for the fibre of $E_2$ over any $\gamma \in \Gal_F$.  Again, to understand the extension $E_2$, it suffices to understand these fibres (as $\tilde Z^\vee$-torsors), and to understand the multiplication maps among them.  The steps above yield the following (somewhat) concise description of $E_{2,\gamma}$.

$E_{2,\gamma}$ is the set of homomorphisms $\chi \From D_{Q,n} \To \CC^\times$ such that
\begin{itemize}
\item
$\chi$ is trivial on the image of $Y_{Q,n}^{\SC}$ via the splitting $s_\psi$.
\item
For every $u \in F^\times$, $\chi(u) = \epsilon( \gamma^{-1} \sqrt[n]{u} / \sqrt[n]{u} )$.
\end{itemize}
Multiplication among fibres is given by usual multiplication, $\chi_1, \chi_2 \mapsto \chi_1 \chi_2$.  The $\tilde Z^\vee$-torsor structure on the fibres is given as follows:  if $\eta \in \tilde Z^\vee$, then 
$$[\eta \ast \chi](d) = \eta(y) \cdot \chi(d) \text{ for all } d \in D_{Q,n} \text{ lying over } y \in Y_{Q,n}.$$

\subsection{Comparison}

Now we describe a map from $\pi_1^{\et} (\gerb{E}_\epsilon(\alg{\tilde G}), \bar z)$ to $E_2$, fibrewise over $\Gal_F$.  From the splitting $s$ (used to define $\bar z$ and restricting to $s_\psi$ on $Y_{Q,n}^{\SC}$), every element of $\bar D_{Q,n}$ can be written uniquely as $s(y) \cdot u$ for some $y \in Y_{Q,n}$ and some $u \in \bar F^\times$.  Such an element $s(y) \cdot u$ is $\Gal_F$-invariant if and only if 
$$\gamma (s(y)) \gamma(u) = s(y) u, \text{ or equivalently } \frac{ \gamma^{-1} u}{u} \cdot \frac{ \gamma^{-1} s}{s}(y) = 1, \text{ for al } \gamma \in \Gal_F.$$
 
Suppose that $\gamma \in \Gal_F$ and $(\tau,1) \in \Hom(\bar z, {}^\gamma \bar z)$.  Define $\chi \From D_{Q,n} \To \mu_n(\CC)$ by
$$\chi( s(y) \cdot u) = \epsilon \left( \gamma^{-1} \sqrt[n]{u} / \sqrt[n]{u} \cdot \tau(y) \right).$$
This makes sense, because $\Gal_F$-invariance of $s(y) \cdot u$ implies
$$\left( \frac{ \gamma^{-1} \sqrt[n]{u}}{\sqrt[n]{u}} \cdot \tau(y) \right)^n = \frac{ \gamma^{-1} u}{u} \cdot \frac{ \gamma^{-1} s}{s}(y) = 1.$$
To see that $\chi \in E_{2,\gamma}$, observe that 
\begin{itemize}
\item
$\chi$ is a homomorphism (a straightforward computation).
\item
If $y \in Y_{Q,n}^{\SC}$ then $\chi(s(y)) = \tau(y) = 1$ since $\tau \in \hat Z$.  
\item
If $u \in F^\times$ then $\chi(u) =  \epsilon( \gamma^{-1} \sqrt[n]{u} / \sqrt[n]{u} )$ by definition.
\end{itemize}

\begin{lemma}
The map sending $(\tau,1)$ to $\chi$, described above, extends uniquely to an isomorphism of $\tilde Z^\vee$-torsors from $\Hom(\bar z, {}^\gamma \bar z)$ to $E_{2,\gamma}$.
\end{lemma}
\proof
If this map extends to an isomorphism of $\tilde Z^\vee$-torsors as claimed, the map must send an element $(\tau, \zeta) \in \Hom(\bar z, {}^\gamma \bar z)$ to the element $\zeta \ast \chi \in E_{2,\gamma}$.  To demonstrate that the map extends to an isomorphism of $\tilde Z^\vee$-torsors, it must only be checked that
$$(\xi \cdot \tau, 1) \text{ and } (\tau, \epsilon(\xi))$$
map to the same element of $E_{2,\gamma}$, for all $\xi \in \hat Z_{[n]}$.  For this, we observe that $(\xi \cdot \tau, 1)$ maps to the character $\chi'$ given by
$$\chi'( s(y) \cdot u) = \epsilon \left( \frac{ \gamma^{-1} \sqrt[n]{u}}{\sqrt[n]{u}} \xi(y) \tau(y) \right) = \epsilon(\xi(y)) \cdot \epsilon \left( \frac{ \gamma^{-1} \sqrt[n]{u}}{\sqrt[n]{u}} \tau(y) \right) = \epsilon(\xi(y)) \cdot \chi(s(y) \cdot u).$$
Thus $\chi' = \epsilon(\xi) \ast \chi$ and this demonstrates the lemma.
\qed

From this lemma, we have a well-defined ``comparison'' isomorphism of $\tilde Z^\vee$-torsors,
$$C_\gamma \From \Hom(\bar z, {}^\gamma z) \To E_{2,\gamma},$$
\begin{equation}
\tag{Comparison}
C_\gamma(\tau, \zeta) (s(y) \cdot u) = \epsilon \left( \frac{ \gamma^{-1} \sqrt[n]{u}}{\sqrt[n]{u}} \cdot \tau(y) \right) \cdot \zeta(y).
\end{equation}

Checking compatibility with multiplication yields the following.
\begin{lemma}
The isomorphisms $C_\gamma$ are compatible with the multiplication maps, yielding an isomorphism of extensions of $\Gal_F$ by $\tilde Z^\vee$,
$$C = C_{\bar z} \From \pi_1^{\et}(\gerb{E}_\epsilon(\alg{\tilde G}), \bar z) \To E_2.$$
\end{lemma}
\proof
Suppose that $(\tau_1, \zeta_1) \in \Hom(\bar z, {}^{\gamma_1} \bar z)$ and $(\tau_2, \zeta_2) \in \Hom(\bar z, {}^{\gamma_2} \bar z)$.  Their product in $\pi_1^{\et}(\gerb{E}_\epsilon(\alg{\tilde G}), \bar z)$ is $(\gamma_2^{-1}(\tau_1) \tau_2, \zeta_1 \zeta_2)$.  We compute
\begin{align*}
C_{\gamma_1 \gamma_2} (\tau_1 \gamma^{-1}(\tau_2), \zeta_1 \zeta_2) (s(y) \cdot u) &= \epsilon \left( \frac{ (\gamma_1 \gamma_2)^{-1}\sqrt[n]{u}}{\sqrt[n]{u}} \cdot  \gamma_2^{-1}(\tau_1(y)) \tau_2(y) \right) \cdot \zeta_1(y) \zeta_2(y) \\
&= \epsilon \left( \frac{ \gamma_2^{-1} \gamma_1^{-1} \sqrt[n]{u}}{\gamma_2^{-1} \sqrt[n]{u}}  \cdot \frac{ \gamma_2^{-1} \sqrt[n]{u}}{\sqrt[n]{u}} \cdot  \gamma_2^{-1}(\tau_1(y)) \tau_2(y) \right)  \\
&\phantom{=} \cdot \zeta_1(y) \zeta_2(y) \\ 
&= \epsilon \left( \gamma_2^{-1} \left( \frac{ \gamma_1^{-1} \sqrt[n]{u}}{\sqrt[n]{u}} \cdot \tau_1(y) \right) \cdot \frac{ \gamma_2^{-1} \sqrt[n]{u}}{\sqrt[n]{u}} \cdot  \tau_2(y) \right)  \\
&\phantom{=} \cdot \zeta_1(y) \zeta_2(y) \\
&= \epsilon \left(  \frac{ \gamma_1^{-1} \sqrt[n]{u}}{\sqrt[n]{u}} \cdot \tau_1(y) \right) \zeta_1(y) \\
&\phantom{=} \cdot \epsilon \left( \frac{ \gamma_2^{-1} \sqrt[n]{u}}{\sqrt[n]{u}} \cdot  \tau_2(y) \right) \zeta_2(y) \\
&= C_{\gamma_1}(\tau_1, \zeta_1)(s(y) \cdot u) \cdot C_{\gamma_2}(\tau_2, \zeta_2)(s(y) \cdot u)
\end{align*}
In the middle step, we use the fact that $\left( \frac{ \gamma_1^{-1} \sqrt[n]{u}}{\sqrt[n]{u}} \cdot \tau_1(y) \right)$ is an element of $\mu_n(F)$, and hence is $\Gal_F$-invariant.  This computation demonstrates compatibility of the isomorphisms $C_\gamma$ with multiplication maps, and hence the lemma is proven.
\qed

\subsection{Independence of base point}

Lastly, we demonstrate that the comparison isomorphisms 
$$C_{\bar z} \From  \pi_1^{\et}(\gerb{E}_\epsilon(\alg{\tilde G}), \bar z) \To E_2$$
depend naturally on the choice of convenient base point.  With the pinned split group $\alg{G}$ fixed, choosing a convenient base point is the same as choosing a splitting of $\bar D_{Q,n}$ which restricts to $s_\psi$.  

So consider two convenient base points $\bar z_1$ and $\bar z_2$, arising from splittings $s_1, s_2$ of $\bar D_{Q,n}$ which restrict to $s_{\psi}$ on $Y_{Q,n}^{\SC}$.  Any isomorphism $\iota$ from $\bar z_1$ to $\bar z_2$ in the gerbe $\gerb{E}_\epsilon(\alg{\tilde G})$ defines an isomorphism
$$\iota \From \pi_1^{\et}(\gerb{E}_\epsilon(\alg{\tilde G}), \bar z_1) \To \pi_1^{\et}(\gerb{E}_\epsilon(\alg{\tilde G}), \bar z_2).$$
See \cite[Theorem 19.6]{MWBig} for details.  In fact, the isomorphism of fundamental groups above does not depend on the choice of isomorphism from $\bar z_1$ to $\bar z_2$; thus one may define a ``Platonic'' fundamental group
$$\pi_1^{\et}(\gerb{E}_\epsilon(\alg{\tilde G}), \bar F)$$
without reference to an object of the gerbe.

\begin{thm}
For any two convenient base points $\bar z_1, \bar z_2$, and any isomorhpism $\iota \From \bar z_1 \To \bar z_2$, we have $C_{\bar z_2} \circ \iota = C_{\bar z_1}$.  Thus $E_2$ is isomorphic to the fundamental group $\pi_1(\gerb{E}_\epsilon(\alg{\tilde G}), \bar F)$, as defined in \cite[Theorem 19.7, Remark 19.8]{MWBig}.
\end{thm}
\proof
Choose any isomorphism from $\bar z_1 = (\hat T, h_1, j_0)$ to $\bar z_2 = (\hat T, h_2, j_0)$ in the gerbe $\gerb{E}_\epsilon(\alg{\tilde G})$.  Here $h_1(1) = s_1$ and $h_2(1) = s_2$, and $j_0(1) = s_\psi$.  Such an isomorphism $\bar z_1 \xrightarrow{\sim} \bar z_2$ is given by an isomorphism $\iota \From \hat T \To \hat T$ of $\hat T$-torsors satisfying the two conditions
$$h_2 \circ \iota = h_1 \text{ and } j_0 \circ p_\ast \iota = j_0.$$
Such an $\iota$ is determined by the element $b = \iota(1) \in \hat T$.  The two conditions above are equivalent to the two conditions
$$b^n = s_1 / s_2 \text{ and } b \in \hat Z.$$

The isomorphism $\bar z_1 \xrightarrow{\sim} \bar z_2$ determined by such a $b \in \hat T$ yields an isomorphism ${}^\gamma \iota \From {}^\gamma \bar z_1 \To {}^\gamma \bar z_2$, for any $\gamma \in \Gal_F$.  The isomorphim ${}^\gamma \iota$ is given by the isomorphism of $\hat T$-torsors from ${}^\gamma \hat T$ to ${}^\gamma \hat T$, which sends $1$ to $\gamma(b)$.

This allows us to describe the isomorphism
$$\iota \From \pi_1^{\et}(\gerb{E}_\epsilon(\alg{\tilde G}), \bar z_1) \To \pi_1^{\et}(\gerb{E}_\epsilon(\alg{\tilde G}), \bar z_2)$$
fibrewise over $\Gal_F$.  Namely, for any $\gamma \in \Gal_F$, and any $f \in \Hom(\bar z_1, {}^\gamma \bar z_1)$, we find a unique element $\iota(f) \in \Hom(\bar z_2, {}^\gamma \bar z_2)$ which makes the following diagram commute.
$$\begin{tikzcd}
\bar z_1 \arrow{r}{f} \arrow{d}{\iota} & {}^\gamma \bar z_1 \arrow{d}{ {}^\gamma \iota} \\
\bar z_2 \arrow{r}{\iota(f)} & {}^\gamma \bar z_2
\end{tikzcd}$$

If $f = (\tau, 1)$, then $\iota(f) = (\tau b / \gamma^{-1} b, 1)$.  Indeed, when $\tau^n = \gamma^{-1} s_1 / s_1$, we have
$$\left( \frac{ \tau b}{\gamma^{-1} b} \right)^n = \frac{ \gamma^{-1} s_1}{s_1} \frac{b^n}{\gamma^{-1} b^n} = \frac{ \gamma^{-1} s_1}{s_1} \frac{s_1}{s_2} \frac{ \gamma^{-1} s_2}{\gamma^{-1} s_1} = \frac{ \gamma^{-1} s_2}{s_2}.$$
Thus $\iota(f) \in \Hom(\bar z_2, {}^\gamma \bar z_2)$ as required.  In this way,
$$\iota \From \pi_1^{\et}(\gerb{E}_\epsilon(\alg{\tilde G}), \bar z_1) \To \pi_1^{\et}(\gerb{E}_\epsilon(\alg{\tilde G}), \bar z_2),$$
is given concretely on each fibre over $\gamma \in \Gal_F$ by
$$\iota(\tau, \zeta) = \left( \tau \cdot \frac{b}{\gamma^{-1} b}, \zeta \right).$$

Note that the conditions $b^n = s_1 / s_2$ and $b \in \hat Z$ uniquely determine $b$ up to multiplication by $\hat Z_{[n]}$.  Since $\hat Z_{[n]}$ is a trivial $\Gal_F$-module, the isomorphism $\iota$ of fundamental groups is independent of $b$.  Finally, we compute, for any $y \in Y_{Q,n}, u \in \bar F^\times$ such that $s_1(y) \cdot u \in D_{Q,n}$, and any $(\tau, \zeta) \in \Hom(\bar z_1, {}^\gamma \bar z_1)$,
\begin{align*}
[C_{\bar z_2} \circ \iota](\tau, \zeta)(s_1(y) \cdot u) &= C_{\bar z_2}( \tau b / \gamma^{-1} b, \zeta)(s_1(y) \cdot u) \\
&= C_{\bar z_2}( \tau \gamma(b) / b, \zeta)(s_2(y) \cdot b^n(y) u) \\
&= \epsilon \left( \frac{ \gamma^{-1} \sqrt[n]{b^n(y) u} }{\sqrt[n]{b^n(y) u}} \cdot \tau(y) \cdot \frac{b(y)}{ \gamma^{-1}(b(y))} \right) \cdot \zeta(y) \\
&= \epsilon \left( \frac{ \gamma^{-1} \sqrt[n]{u}}{\sqrt[n]{u}} \cdot \tau(y) \right) \cdot \zeta(u) \\
&= C_{\bar z_1}(\tau, \zeta)(s_1(y) \cdot u).
\end{align*}
\qed

As noted in the introduction, this demonstrates compatibility between two approaches to the L-group.
\begin{corollary}
The L-group defined in \cite{MWBig} is isomorphic to the L-group defined in \cite{GanGao}, for all pinned split reductive groups over local or global fields.
\end{corollary}

\printbibliography

\end{document}